\documentclass[submit]{smj}

\usepackage{amsthm}
\usepackage{amsmath}
\usepackage{natbib}
\usepackage{graphicx}
\usepackage{amssymb}
\usepackage[table]{xcolor}
\usepackage{color}
\usepackage{amsfonts}

\usepackage{soul}

\usepackage{enumerate}
\usepackage{stmaryrd}
\usepackage{varioref,exscale,relsize}

\usepackage{placeins}

\usepackage{multirow}
\usepackage{latexsym}
\usepackage{multicol}
\usepackage{dsfont,bbold}

\newcommand{\e}{\varepsilon}

\newcommand{\R}{\mathbb{R}}

\Author{Meili Baragatti \Affil{1}
        Karine Bertin \Affil{2},
        Emilie Lebarbier \Affil{3},
        and Cristian Meza\Affil{2}
} \AuthorRunning{Meili Baragatti \textrm{et al.}}

\Affiliations{

\item SupAgro-INRA UMR MISTEA,
      Montpellier,
      FRANCE

\item CIMFAV-Facultad de Ingenier\'ia,
      Universidad de Valpara\'iso,
      Valpara\'iso,
      Chile

\item AgroParisTech/INRA UMR518,
      Paris 5e,
      France

}   

\CorrAddress{Karine Bertin,
             CIMFAV-Facultad de Ingenier\'ia,
             Universidad de Valpara\'iso,
             Valpara\'iso,
             Chile}
\CorrEmail{karine.bertin@uv.cl} \CorrPhone{(+56)\;32\;260\;36\;23}
\CorrFax{(+56)\;32\;260\;36\;18}

\Title{A Bayesian approach for the segmentation of series corrupted by a functional part}
\TitleRunning{A Bayesian approach for the segmentation of series corrupted by a functional part}

\Abstract{We propose a Bayesian approach to detect multiple change-points in a piecewise-constant signal corrupted by a functional part corresponding to  environmental or experimental disturbances. The piecewise constant part (also called segmentation part) is expressed as the product of a lower triangular matrix by a sparse vector. The functional part is a linear combination of functions from a large dictionary.  A Stochastic Search Variable Selection approach is used to obtain sparse estimations of the segmentation parameters (the change-points and the means over the segments) and of the functional part. The performance of our proposed method is assessed using simulation experiments. Applications to two real datasets from geodesy and economy fields are also presented.}


\begin{document}

\maketitle


\section{Introduction}

The problem of detecting multiple change-points in signals arises in many fields such as biology \citep{BH:2004}, geodesy (\citealp{W2003}; \citealp{BCLM:2014}), meteorology (\citealp{CM2004}; \citealp{F:2006}; \citealp{WFR:2011}; \citealp{R:2013}) or astronomy \citep{DTS:2007} among others. In addition to change-points, we may observe environmental or experimental disturbances (for instance geophysical signals or climatic effects) which need to be taken into account in the models. Since the form of these disturbances is in general unknown or partially unknown, it seems natural to model them as a functional part that has to be estimated. Our goal in this article is to develop a Bayesian approach that allows us to both estimate the segmentation part (the change-points and the means over the segments) and the functional part.
The Bayesian approach has the advantage that expert knowledge can be introduced in the models through prior distributions. This can be useful in multiple change-points problems where change-points can be related to specific events such as instrumental changes, earthquakes, very hot years or months or economic crisis for example. Moreover, posterior distributions allow us for a quantification of the uncertainty, giving in particular posterior probabilities or credible intervals for the positions of change-points or the functional part. This is of particular interest for practitioners.

Several methods have been proposed in a Bayesian framework for the multiple change-points problem. These methods are based, mostly, on reversible jump Markov Chain Monte Carlo algorithms (\citealp{LL2000}; \citealp{BH:2004}; \citealp{Tai:2010}),  Stochastic search Variable Selection \citep{DTS:2007},  dynamic programming recursions \citep{R:2013} or non-parametric Bayesian approaches (\citealp{Martinez:14} and references there in). All these Bayesian methods deal with the multiple change-points problem but they do not consider the presence of functional disturbances.
However, as illustrated in \cite{PLBR11} and \cite{BCLM:2014} in simulation and real examples, taking into account the functional part in the segmentation model can be crucial for an accurate change-point detection and interesting information can be extracted from the form of the functional part.

We propose a novel Bayesian method to detect multiple change-points in a piecewise-constant signal corrupted by a functional part, where
 the functional part is estimated using a dictionary approach \citep{tsybakov} and the segmentation part is treated as a sparse problem. More precisely, concerning the segmentation, we follow \cite{HL:2010} by expressing the piecewise constant part of the model as a product of a lower triangular matrix by a sparse vector (which non-zero coordinates correspond to change-points positions). In addition, the functional part is represented as a linear combination of functions from a dictionary. Since a large variety of functions can be included in the dictionary, this leads generally to a sparse representation of the functional part in terms of functions from the dictionary. Hence, a Stochastic search Variable Selection approach can be used to estimate the sparse vectors, that is, both the location of the change-points and the functional part \citep{GMC1997}.

In the simulation studies and for real datasets, we obtain good results for both the segmentation and the functional parts. On a GPS series from an Australian station and on a series of daily records of Mexican Peso/US Dollar exchange rate, expected change-points are recovered. Moreover the benefits of the Bayesian approach are also illustrated in the simulation and real data studies. In particular, for the GPS series, the use of prior knowledge about instrumental changes enables us to detect some relevant change-points and, in the functional part, periodic components suggested in previous works.

The remainder of the paper is organized as follows. Section \ref{sec:model} presents the hierarchical Bayesian model considered, Section \ref{sec:MCMC} outlines the procedure used to estimate the model parameters. In Section \ref{sec:simu}, the performance of the proposed method is studied through simulations. The procedure is illustrated on the two real datasets in Section \ref{sec:appli}, and finally Section \ref{sec:discussion} discusses it.

\section{Model}\label{sec:model}

\subsection{Segmentation model with functional part}\label{modelsemi}

We observe a series $\mathbf{Y}=(Y_1,\ldots,Y_n)'$ that satisfies
\begin{equation}\label{model1}
Y_t=\mu_k+f(x_t)+\e_t, \quad \forall t\in I_k=(\tau_{k-1},\tau_k], k\in\{1,\ldots,K\},
\end{equation}
where $K$ is the total number of segments of the  series is unknown, the $\e_t$ are i.i.d centered Gaussian variables with variance $\sigma^2$, $x_t$ is a covariate (the simple one is the time $t$), $f$ is an unknown function to be estimated, $\tau_k$ is the $k$th change-point, $\mu_{k}$ is the mean of the series on the segment $I_k$. We use the convention $\tau_0=0$ and $\tau_K=n$.

A classical approach in non-parametric framework is to expand the functional part $f$ with respect to orthonormal basis, such as Fourier or wavelet ones (see \citealp{picard} and references therein). Following \cite{tsybakov} or \cite{BCLM:2014}, we choose here to adopt a dictionary approach,  that consists in finding an over-complete representation of $f$. More precisely, we expand $f$ with respect to a large family of functions $(\phi_j)_{j=1,\ldots,M}$, named dictionary, that can  for example be the union of two orthonormal basis. Then $f$ is assumed to be of the form
\begin{equation*}\label{model2}
f(x)=\sum_{j=1}^M \lambda_j\phi_j(x),
\end{equation*}
where $\boldsymbol{\lambda}=(\lambda_1,\ldots,\lambda_M)'\in\R^M$ is a vector of coordinates of $f$ in the dictionary and
\begin{equation*}
(f(x_1),\ldots,f(x_n))'=\mathbf{F}\boldsymbol{\lambda},
\end{equation*}
where $\mathbf{F}$ is the $n\times M$ matrix $\mathbf{F}=(\phi_j(x_i))_{i,j}$.
Note that since large dictionaries are considered, this allows us to obtain a sparse representation of the function $f$, that is the vector $\boldsymbol{\lambda}$ is expected to be with few non-zero coordinates.

To estimate the change-points in the series, we follow the strategy proposed by \cite{HL:2010}, which consists in reframing this task in a variable selection context. We denote by $\mathbf{X}$ the $n\times n$ lower triangular matrix having only $1$'s on the diagonal and below it. We consider the $n\times 1$ vector $\boldsymbol{\beta}$ with only $K$ non-zero coefficients at positions $(\tau_k+1)_{k=0,\ldots,K-1}$ with $\beta_{\tau_k+1}=\mu_{k+1}-\mu_k$ and using the convention $\mu_0=0$. Note that the segmentation (the change-points $\tau_k$ and the means $\mu_k$) will be recovered by the vector $\boldsymbol{\beta}$.

The model (\ref{model1}) can then be rewritten as follows
\begin{equation*}\label{model3}
\mathbf{Y}=\mathbf{X}\boldsymbol{\beta}+\mathbf{F}\boldsymbol{\lambda}+\boldsymbol{\e},
\end{equation*}
where $\boldsymbol{\e}=(\e_1,\ldots,\e_n)'$. Our objective is  now to estimate the parameters $\boldsymbol{\beta}$, $\boldsymbol{\lambda}$ and $\sigma^2$. Since both $\boldsymbol{\beta}$ and $\boldsymbol{\lambda}$ vectors are expected to be sparse, we propose to use Bayesian methods of variable selection for their estimation.

\subsection{Bayesian hiercharchical framework}\label{subsec:prior}

Following \cite{GeorgeMcCulloch}, we first introduce latent variables $\boldsymbol{\gamma}$ and $\mathbf{r}$ to identify non-null components of the vectors $\boldsymbol{\beta}$ and $\boldsymbol{\lambda}$.
The vector $\boldsymbol{\gamma}=(\gamma_1,\ldots,\gamma_n)$ is such that $\gamma_i=I_{\{\beta_i\neq 0\}}$, where $I$ denotes the indicator function and the vector $\mathbf{r}=(r_1,\ldots,r_M)$ satisfies $r_j=I_{\{\lambda_j\neq 0\}}$. The number of non-zero coordinates of $\boldsymbol{\gamma}$ and $\mathbf{r}$ are $d_{\boldsymbol{\gamma}}=K$ and $d_{\mathbf{r}}$ respectively.
The product $\mathbf{X}\boldsymbol{\beta}$ is equal to $\mathbf{X}_{\boldsymbol{\gamma}}\boldsymbol{\beta}_{\boldsymbol{\gamma}}$ where $\mathbf{X}_{\boldsymbol{\gamma}}$ is the $n\times d_{\boldsymbol{\gamma}}$ matrix containing only the $j$ columns of $X$ such that $\gamma_j$ is non-zero and $\boldsymbol{\beta}_{\boldsymbol{\gamma}}$ is a $d_{\boldsymbol{\gamma}}\times 1$ vector containing only the non-zero coefficients of $\boldsymbol{\beta}$. Similarly, we can express $\mathbf{F}\boldsymbol{\lambda}$ as $\mathbf{F}_{\mathbf{r}}\boldsymbol{\lambda}_{\mathbf{r}}$ where $\mathbf{F}_{\mathbf{r}}$ is a $n\times d_{\mathbf{r}}$ matrix and $\boldsymbol{\lambda}_{\mathbf{r}}$ a $d_{\mathbf{r}}\times 1$ vector. The model (\ref{model1}) can be then rewritten as

\begin{equation*}\label{model4}
\mathbf{Y}=\mathbf{X}_{\boldsymbol{\gamma}}\boldsymbol{\beta}_{\boldsymbol{\gamma}}+\mathbf{F}_{\mathbf{r}}\boldsymbol{\lambda}_{\mathbf{r}}+\boldsymbol{\e}
\end{equation*}
where the parameters to estimate are $\boldsymbol{\theta}=\{\boldsymbol{\beta}_{\boldsymbol{\gamma}},\boldsymbol{\gamma},\boldsymbol{\lambda}_{\mathbf{r}},r,\sigma^2\}$.

Then, as usual in a Bayesian context, these parameters are treated as random variables, assumed here to be independent, and we consider the following prior distributions. The $\gamma_i$ are independent Bernoulli variables with parameter $0\le\pi_i\le 1$ for $i=2,\dots,n$ and with $\pi_1=1$ by convention. The $r_j$ are also independent Bernoulli variables with parameter $0\le\eta_j\le 1$  for $j=1,\dots,M$. Then the noise parameter follows a Jeffrey distribution, $\pi(\sigma^2)\propto \sigma^{-2}$. The conditional distribution of $\boldsymbol{\beta}_{\boldsymbol{\gamma}}|\boldsymbol{\gamma},\sigma^2$ is the classical $g$-prior of \cite{ZE1986} given by $\boldsymbol{\beta}_{\boldsymbol{\gamma}}|\boldsymbol{\gamma},\sigma^2\sim\mathcal{N}_{d_{\boldsymbol{\gamma}}}\left(0,c_1\sigma^2\left(\mathbf{X}_{\boldsymbol{\gamma}}'\mathbf{X}_{\boldsymbol{\gamma}}\right)^{-1}\right)$. Finally the conditional distribution of $\boldsymbol{\lambda}_{\mathbf{r}}|\mathbf{r},\sigma^2 $ is also a $g$-prior, with $\boldsymbol{\lambda}_{\mathbf{r}}|\mathbf{r},\sigma^2\sim \mathcal{N}_{d_{\mathbf{r}}}\left(0,c_2\sigma^2\left(\mathbf{F}_{\mathbf{r}}'\mathbf{F}_{\mathbf{r}}\right)^{-1}\right)$.

The posterior distribution of $\boldsymbol{\theta}$ has the following expression
\begin{eqnarray}\label{jointposterior}
\pi(\boldsymbol{\theta}|\mathbf{Y})&=& \frac{\pi(\mathbf{Y}|\boldsymbol{\theta})\pi(\boldsymbol{\beta}_{\boldsymbol{\gamma}}|\boldsymbol{\gamma},\sigma^2)\pi(\boldsymbol{\lambda}_{\mathbf{r}}|\mathbf{r},\sigma^2)\pi(\boldsymbol{\gamma})\pi(\mathbf{r})\pi(\sigma^2)}{\pi(\mathbf{Y})},
\end{eqnarray}
where
\begin{equation*}\pi(\mathbf{Y}|\boldsymbol{\theta})=\left(\frac{1}{2\pi\sigma^2} \right)^{\frac{n}{2}}\exp\left(-\frac{1}{2\sigma^2} \left(\mathbf{Y}-\mathbf{X}_{\boldsymbol{\gamma}}\boldsymbol{\beta}_{\boldsymbol{\gamma}} -\mathbf{F}_{\mathbf{r}}\boldsymbol{\lambda}_{\mathbf{r}}\right)^\prime \left(\mathbf{Y}-\mathbf{X}_{\boldsymbol{\gamma}}\boldsymbol{\beta}_{\boldsymbol{\gamma}} -\mathbf{F}_{\mathbf{r}}\boldsymbol{\lambda}_{\mathbf{r}}\right) \right).
\end{equation*}

\section{MCMC Schemes}\label{sec:MCMC}

A classical approach for the computational scheme would be to estimate the whole parameters at the same time $(\boldsymbol{\beta}_{\boldsymbol{\gamma}},\boldsymbol{\gamma},\boldsymbol{\lambda}_{\mathbf{r}},\mathbf{r},\sigma^2)$ using a Metropolis-within-Gibbs algorithm combined with the grouping (or blocking) technique of \citet{Liu}. Indeed, $\boldsymbol{\beta}_{\boldsymbol{\gamma}}$ and $\boldsymbol{\gamma}$, as well as $\boldsymbol{\lambda}_{\mathbf{r}}$ and $\mathbf{r}$, cannot be considered separately. An iteration of the algorithm would be made of three steps: update of $\boldsymbol{\beta}_{\boldsymbol{\gamma}},\boldsymbol{\gamma} \mid \boldsymbol{\lambda}_{\mathbf{r}},\mathbf{r},\sigma^2, \mathbf{Y}$, update of $\boldsymbol{\lambda}_{\mathbf{r}},\mathbf{r} \mid \boldsymbol{\beta}_{\boldsymbol{\gamma}},\boldsymbol{\gamma}, \sigma^2, \mathbf{Y}$ and update of $\sigma^2 \mid \boldsymbol{\beta}_{\boldsymbol{\gamma}},\boldsymbol{\gamma},\boldsymbol{\lambda}_{\mathbf{r}},\mathbf{r}, \mathbf{Y}$.  However some drawbacks are associated with this algorithm. Firstly the update rates for $(\boldsymbol{\beta}_{\boldsymbol{\gamma}},\boldsymbol{\gamma})$ and $(\boldsymbol{\lambda}_{\mathbf{r}},\mathbf{r})$ will be very low since it is difficult to make good proposals for both $\boldsymbol{\beta}_{\boldsymbol{\gamma}}$ and $\boldsymbol{\gamma}$ or for both $\boldsymbol{\lambda}_{\mathbf{r}}$ and $\mathbf{r}$. Secondly, the post-hoc interpretation of the vectors $\boldsymbol{\beta}_{\boldsymbol{\gamma}}$ and $\boldsymbol{\lambda}_{\mathbf{r}}$ obtained at each iteration of this algorithm are not relevant. Indeed, they are not associated with the same change-points and functions from one iteration to the next.
Besides, as explained in \cite{LL2000} (see Section 4), in a Bayesian segmentation framework, the posterior mean of $\boldsymbol{\beta}_{\boldsymbol{\gamma}}$, obtained from the posterior distribution of $(\boldsymbol{\beta}_{\boldsymbol{\gamma}},\boldsymbol{\gamma})$, is uninteresting. The interpretation of this mean is not obvious since it is calculated over all the possible configurations of change-points. A solution is to use the posterior distribution of $\boldsymbol{\beta}_{\boldsymbol{\gamma}}$ conditionally to $\boldsymbol{\gamma}$. We have the same drawback for the functional part.

For these two reasons, we propose the following two-step strategy: the first step aims at detecting the positions of the change-points and at selecting the functions, that is, to estimate the latent vectors $\boldsymbol{\gamma}$ and $\mathbf{r}$. To this end, the parameters $\boldsymbol{\beta}_{\boldsymbol{\gamma}}$, $\boldsymbol{\lambda}_{\mathbf{r}}$ and $\sigma^2$ can be considered as nuisance parameters, and we use the joint posterior distribution integrated with respect to $\boldsymbol{\beta}_{\boldsymbol{\gamma}}$, $\boldsymbol{\lambda}_{\mathbf{r}}$ and $\sigma^2$. This can be viewed as a collapsing technique, see \citet{Liu} and \citet{vanDyk}. In the second part, we estimate $\boldsymbol{\beta}_{\boldsymbol{\gamma}}$, $\boldsymbol{\lambda}_{\mathbf{r}}$ and $\sigma^2$, conditionally to $\boldsymbol{\gamma}$ and $\mathbf{r}$. The MCMC scheme would then be as follows:

\begin{enumerate}
\item Estimation of $\boldsymbol{\gamma}$ and $\mathbf{r}$: use of a Metropolis-Hastings algorithm to draw from the joint posterior distribution $\pi(\boldsymbol{\gamma},\mathbf{r}|\mathbf{Y})$ integrated with respect to $\boldsymbol{\beta}_{\boldsymbol{\gamma}}$, $\boldsymbol{\lambda}_{\mathbf{r}}$ and $\sigma^2$.
\item Estimation of $\boldsymbol{\beta}_{\boldsymbol{\gamma}}, \boldsymbol{\lambda}_{\mathbf{r}}$ and $\sigma$: given the estimates $\hat{\boldsymbol{\gamma}}$ and $\hat{\mathbf{r}}$, use of a Gibbs sampler algorithm.
\end{enumerate}
In the following subsections, we give some details of both steps.

\subsection{Metropolis-Hastings algorithm}\label{subsec:MH}
The joint posterior distribution integrated with respect to $\boldsymbol{\beta}_{\boldsymbol{\gamma}}$, $\boldsymbol{\lambda}_{\mathbf{r}}$ and $\sigma^2$ is the following (see details in Appendix \ref{App:AppendixA}):
\begin{equation*}\label{integratedposterior}
\pi(\boldsymbol{\gamma},\mathbf{r}|\mathbf{Y})\propto (1+c_1)^{-d_{\boldsymbol{\gamma}}/2}\pi(\boldsymbol{\gamma})\pi(\mathbf{r}) g(\boldsymbol{\gamma},\mathbf{r},\mathbf{Y}),
\end{equation*}
where
\begin{align*}
g(\boldsymbol{\gamma},\mathbf{r},\mathbf{Y})=&\left(\frac{\left| \left( \mathbf{F}_{\mathbf{r}}^\prime \left(\mathbf{U}_{\boldsymbol{\gamma}}^{-1}+\frac{I}{c_2} \right) \mathbf{F}_{\mathbf{r}}\right)^{-1}  \right|}{|c_2(\mathbf{F}_{\mathbf{r}}^\prime \mathbf{F}_{\mathbf{r}})^{-1}|} \right)^{1/2}\\
&\times \left[ \frac{1}{2} \mathbf{Y}^\prime\left( \mathbf{U}_{\boldsymbol{\gamma}}^{-1}-\mathbf{U}_{\boldsymbol{\gamma}}^{-1} \mathbf{F}_{\mathbf{r}}\left(\mathbf{F}_{\mathbf{r}}^\prime \left(\mathbf{U}_{\boldsymbol{\gamma}}^{-1}+\frac{I}{c_2} \right) \mathbf{F}_{\mathbf{r}} \right)^{-1}\mathbf{F}_{\mathbf{r}}^\prime \mathbf{U}_{\boldsymbol{\gamma}}^{-1} \right) \mathbf{Y}  \right]^{-n/2},
\end{align*}
and
\begin{equation*}
\mathbf{U}_{\boldsymbol{\gamma}}=\left(I-\frac{c_1}{1+c_1} \mathbf{X}_{\boldsymbol{\gamma}}\left(\mathbf{X}_{\boldsymbol{\gamma}}^\prime \mathbf{X}_{\boldsymbol{\gamma}}\right)^{-1} \mathbf{X}_{\boldsymbol{\gamma}}^\prime \right)^{-1}.
\end{equation*}

To sample from $\pi(\boldsymbol{\gamma},\mathbf{r}|\mathbf{Y})$, a Metropolis-Hastings algorithm is used.
At iteration $t$, a candidate $(\boldsymbol{\gamma}^*,\mathbf{r}^*)$ is proposed from $(\boldsymbol{\gamma}^{(t)},\mathbf{r}^{(t)})$, and using symmetric transition kernel the acceptance rate is:
\begin{equation*}
\rho\big((\boldsymbol{\gamma}^{(t)},\mathbf{r}^{(t)});(\boldsymbol{\gamma}^*,\mathbf{r}^*)\big)  =
\min \Bigg\{1,\frac{\pi(\boldsymbol{\gamma}^*,\mathbf{r}^*|\mathbf{Y})}{\pi(\boldsymbol{\gamma}^{(t)},\mathbf{r}^{(t)}|\mathbf{Y})} \Bigg\}.
\end{equation*}

To have a symmetric transition kernel, two kinds of proposals are used (each one of probability 1/2): either $k$ components of $\boldsymbol{\gamma}^{(t)}$ are randomly changed, or $l$ components of $\mathbf{r}^{(t)}$ are randomly changed (a value of 0 is switched to 1, and conversely). We modify only one of the two latent vectors at each iteration since proposals with modifications of both vectors have too low acceptance rates. In this algorithm, by convention, we suppose $\gamma_1=1$ (time 0 is a change-point, corresponding to $\tau_0=0$) and $r_1=1$ (the constant function is always selected).

The number of iterations of this algorithm is $b+m$, where $b$ corresponds to the burn-in period. Then $\boldsymbol{\gamma}$ and $\mathbf{r}$ are estimated using the sequences $\{\boldsymbol{\gamma}^{(t)}\}$ and $\{\mathbf{r}^{(t)}\}$, for $t=b+1,\ldots,b+m$. The most relevant positions for the change-points and the most relevant functions for the functional part are those which are supported by the data and prior information. In other words, they are those corresponding to the $\boldsymbol{\gamma}$ and $\mathbf{r}$ components with higher posterior probabilities. In practice, the selected components are those such their posterior probability  is higher than a given threshold.
As in \cite{MPRR04} or \cite{MPR06}, we choose a threshold that minimized a loss function. Here the considered loss function is the sum of the False Discovery and of the False Negative ($FD+FN$), leading to a threshold of $1/2$. This also corresponds to the selection of the median probability model in \cite{BB2004}, which has been shown to have greater predictive power than the most probable model, under many circumstances.

\subsection{Gibbs sampler algorithm}
Once $\boldsymbol{\gamma}$ and $\mathbf{r}$ have been estimated, our goal is to estimate $\boldsymbol{\beta}_{\boldsymbol{\gamma}}$, $\boldsymbol{\lambda}_{\mathbf{r}}$ and $\sigma^2$ from the distribution
 $\pi(\boldsymbol{\beta}_{\boldsymbol{\gamma}},\boldsymbol{\lambda}_{\mathbf{r}},\sigma^2|\mathbf{r},\boldsymbol{\gamma},\mathbf{Y})\propto \pi(\boldsymbol{\beta}_{\boldsymbol{\gamma}},\boldsymbol{\lambda}_{\mathbf{r}},\sigma^2,\mathbf{r},\boldsymbol{\gamma}|\mathbf{Y}).$

A Gibbs sampler algorithm is then used. At each iteration, the three parameters should be drawn from its full conditional distribution given by:
\begin{align*}
 &\boldsymbol{\beta}_{\boldsymbol{\gamma}}|\boldsymbol{\lambda}_{\mathbf{r}},\sigma^2,\mathbf{r},\boldsymbol{\gamma},\mathbf{Y}\sim \mathcal{N}_{d_{\boldsymbol{\gamma}}}\left(\frac{\mathbf{T}_{\boldsymbol{\gamma}} \mathbf{X}_{\boldsymbol{\gamma}}^\prime(\mathbf{Y}-\mathbf{F}_{\mathbf{r}}\boldsymbol{\lambda}_{\mathbf{r}})}{\sigma^2}, \mathbf{T}_{\boldsymbol{\gamma}} \right),\\
 & \boldsymbol{\lambda}_{\mathbf{r}}|\boldsymbol{\beta}_{\boldsymbol{\gamma}},\sigma^2,\mathbf{r},\boldsymbol{\gamma},\mathbf{Y}\sim \mathcal{N}_{d_{\mathbf{r}}}\left( \frac{\mathbf{W}_{\mathbf{r}}\mathbf{F}_{\mathbf{r}}^\prime(\mathbf{Y}-\mathbf{X}_{\boldsymbol{\gamma}}\boldsymbol{\beta}_{\boldsymbol{\gamma}})}{\sigma^2}, \mathbf{W}_{\mathbf{r}}\right),\\
 &\sigma^2|\boldsymbol{\beta}_{\boldsymbol{\gamma}},\boldsymbol{\lambda}_{\mathbf{r}},\mathbf{r},\boldsymbol{\gamma},\mathbf{Y}\sim IG\left(a,\frac{b}{2}\right),
\end{align*}
where $\mathbf{T}_{\boldsymbol{\gamma}} =\sigma^2\left[\frac{1+c_1}{c_1} \mathbf{X}_{\boldsymbol{\gamma}}^\prime \mathbf{X}_{\boldsymbol{\gamma}} \right]^{-1}$, $\mathbf{W}_{\mathbf{r}}=\sigma^2\left[\frac{1+c_2}{c_2} \mathbf{F}_{\mathbf{r}}^\prime \mathbf{F}_{\mathbf{r}} \right]^{-1}$, $a=\frac{n}{2}+\frac{d_{\boldsymbol{\gamma}}}{2}+\frac{d_{\mathbf{r}}}{2}$ and
\begin{equation*}
 b=\left(\mathbf{Y}-\mathbf{X}_{\boldsymbol{\gamma}}\boldsymbol{\beta}_{\boldsymbol{\gamma}}-\mathbf{F}_{\mathbf{r}}\boldsymbol{\lambda}_{\mathbf{r}}\right)^\prime \left(\mathbf{Y}-\mathbf{X}_{\boldsymbol{\gamma}}\boldsymbol{\beta}_{\boldsymbol{\gamma}}-\mathbf{F}_{\mathbf{r}}\boldsymbol{\lambda}_{\mathbf{r}}\right) + \boldsymbol{\beta}_{\boldsymbol{\gamma}}^\prime \left(\frac{\mathbf{X}_{\boldsymbol{\gamma}}^\prime \mathbf{X}_{\boldsymbol{\gamma}}}{c_1}\right)\boldsymbol{\beta}_{\boldsymbol{\gamma}} + \boldsymbol{\lambda}_{\mathbf{r}}^\prime \left(\frac{\mathbf{F}_{\mathbf{r}}^\prime \mathbf{F}_{\mathbf{r}}}{c_2}\right)\boldsymbol{\lambda}_{\mathbf{r}}.
\end{equation*}

To estimate $\boldsymbol{\beta}_{\boldsymbol{\gamma}}$, $\boldsymbol{\lambda}_{\mathbf{r}}$ and $\sigma^2$, empirical posterior means are computed using only post-burn-in iterations. Using the estimators $\hat{\boldsymbol{\beta}}$ and $\hat{\boldsymbol{\lambda}}$  of $\boldsymbol{\boldsymbol{\beta}}$ and $\boldsymbol{\lambda}$, we then obtain the estimator $\hat{f}(\cdot)=\sum_{j=1}^M\hat{\lambda}_j\phi_j(\cdot)$ of the function $f$,
the estimators $\hat{\tau}_k$ of the change-points and the estimators $\hat{\mu}_k$ of the means (see Section~\ref{modelsemi}). The estimated number of change-points is given by $\hat{K}=\sum_{i=1}^n I_{\hat{\beta}_i\neq 0}$.

\section{Simulation study}\label{sec:simu}

In this section, we conduct a simulation study to assess the performance of our proposed method to estimate both the "parametric" part (the segmentation part) and the``non-parametric" part (the functional part). Our method is called {\it{SegBayes\_SP}} for``semi-parametric". Moreover in order to investigate the impact of taking into account the functional part on the
estimation of the segmentation part, we compare the results of {\it{SegBayes\_SP}} with those of the same Bayesian procedure that includes only the segmentation part in the model (i.e. if the model is supposed to be $Y_t = \mu_k +\e_t, \quad \forall t\in I_k=(\tau_{k-1},\tau_k], k\in\{1,\ldots,K\}$). This case is called {\it{SegBayes\_P}} for``parametric". The estimation is still obtained using a Metropolis-Hastings algorithm to estimate $\boldsymbol{\gamma}$ followed by a Gibbs sampler to estimate $\boldsymbol{\beta}_{\boldsymbol{\gamma}}$ and $\sigma^2$. Section~\ref{sec:design} contains our simulation design, the parameters needed for the procedures and the quality criteria. Section~\ref{sec:results} gives the results.

\subsection{Simulation design, parameters of the procedures and quality criteria}\label{sec:design}

\paragraph{Simulation design.}
We consider model (\ref{model1}) with $x_t=t$ and the function $f$ which is a mixture of a sine function with three peaks:
\begin{equation}\label{bias_function}
f(t)=0.3 \times \sin\Big(2\pi \frac{t}{20}\Big) + 1.5 I_{t=0.1\times n} - 2 I_{t=0.5\times n} + 3 I_{t=0.6\times n}.
\end{equation}
Note that this function contains both smooth components and local irregularities (see plot ($h$) of Figure \ref{Exemple_SP_0.1}). We simulate 100 series of length $n=100$ with $K=4$ segments and the mean of each segment takes a value in $\{0,1,2,3,4,5\}$ randomly. The positions of the three change-points are randomly chosen with the following constraints: they are positioned at a distance from the peaks of at least 3, and each segment is at least of length 5. In order to consider several change-point detection difficulties, four levels of noise $\sigma$ are considered: $\sigma \in \{0.1,0.5,1,1.5\}$ (the more  $\sigma$ increases, the more difficult is the detection).

\paragraph{Parameters.}
Several parameters or quantities need to be fixed in the procedures. For the procedure {\it{SegBayes\_SP}}, we consider the following dictionary which contains 151 functions: 128 Haar functions ($t \mapsto 2^{7/2}\mathbb{I}_{[0,1]}\big(\frac{2^7t}{100}-k\Big)$, $k=0,\ldots,2^7-1$), the Fourier functions ($t \mapsto \sin\big(2\pi j \frac{t}{n}\big)$, $t \mapsto \cos\big(2\pi j \frac{t}{n}\big)$, $j=1,\ldots,10$), the functions $t \mapsto t$ and $t \mapsto t^2$, and the constant function. Note that Table \ref{indicesFunctions} gives the indexes of some of these functions (those of $f$ given by (\ref{bias_function}) are in particular 11, 51, 61 and 110).

For both procedures, we use the same prior parameters. With respect to the Metropolis-Hastings algorithm and the Gibbs sampler, we run each one for 20000 iterations including 5000 burn-in iterations. The parameters $c_1$ and $c_2$ are fixed to 50, which is quite standard and recommended for instance by \cite{SmithKohn}. The initial numbers of segments and dictionary functions are 3. The number of change-points proposed to be changed at each iteration is 2, as well as the number of functions of the dictionary.
The initial probability for each position of change-point is 0.01, as well as the initial probability for each function of the dictionary.
To select relevant change-points and functions, we used a threshold equals to 1/2 (see section \ref{subsec:MH} for a justification of this threshold).

\paragraph{Quality criteria.}
To study the performance of the procedures, the same criteria are used for both the segmentation and the functional parts. The first one is the root mean squared error (RMSE) that allows us to assess the quality of the estimation. The other ones (the false discovery rate ($FDR$) and the false negative rate ($FNR$)) evaluate the quality of detection of the change-points and selection of functions. These criteria are detailed in the following:
\begin{itemize}
\item For the segmentation part,
\begin{itemize}
\item the root mean squared distance between the true mean and its estimate is $RMSE(\mu)=\sqrt{\frac{1}{n}\sum_{t=1}^n(\mu(t)-\widehat{\mu}(t))^2}$, with $\mu(t)=\mu_k$ for $t\in I_k=(\tau_{k-1},\tau_k]$, $k=1,\ldots,K$ and $\widehat{\mu}(t)=\widehat{\mu_k}$ for $t\in \widehat{I_k}=(\widehat{\tau}_{k-1},\widehat{\tau}_k]$, $k=1,\ldots,\widehat{K}$.
\item the proportion of erroneously detected change-points among detected change-points, denoted $FDRbp$, and the proportion of undetected change-points among true change-points denoted $FNRbp$.
\end{itemize}
\item For the functional part,
\begin{itemize}
\item the root mean squared distance between $f$ and its estimate:
$RMSE(f)=\sqrt{\frac{1}{n}\sum_{t=1}^n(f(t)-\widehat{f}(t))^2}$.
\item the proportion of erroneously detected functions among detected functions, denoted $FDR(f)$, and the proportion of undetected functions among true functions denoted $FNR(f)$.
\end{itemize}
Note that a perfect segmentation results in both null $FDRbp$ and $FNRbp$, as well as both null $FDR(f)$ and $FNR(f)$ is equivalent to a perfect selection of the functions.
\end{itemize}
The averages of these criteria over the 100 simulations are considered. Note that we also look at the estimation of the standard deviation of the noise $\widehat{\sigma}$.

\subsection{Results}\label{sec:results}
Both procedures have been implemented in R (development core team, 2009) on quad-core processors 2.8GHz Intel X5560 Xeons with 8MB cache size and 32 GB total physical memory. In particular, the standard Metropolis-Hastings algorithm for procedure \textit{SegBayes\_SP} took 169 minutes, for the 400 simulations
(100 series for 4 noises, so it took approximately 25 seconds for each simulation).\\
In the following, we first analyse the results obtained with the 100 simulated series and then we draw some remarks from the study of a particular series.

\begin{figure}[!htp]
    \begin{center}
    \includegraphics[width=13cm]{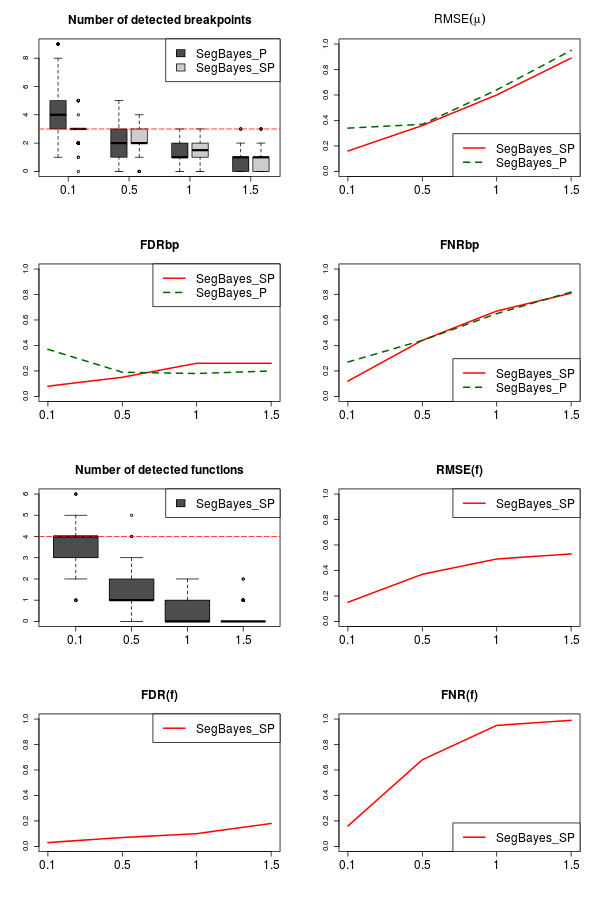}
    \caption{\small{Average quality criteria on 100 simulated series with $\sigma=0.1$, $\sigma=0.5$, $\sigma=1$ and $\sigma=1.5$, for {\it SegBayes\_SP} (semi-parametric model) and {\it SegBayes\_P} (parametric model).}}
    \label{simulations100graph}
    \end{center}
    \end{figure}

\begin{figure}[!htp]
    \begin{center}
    \includegraphics[width=10cm]{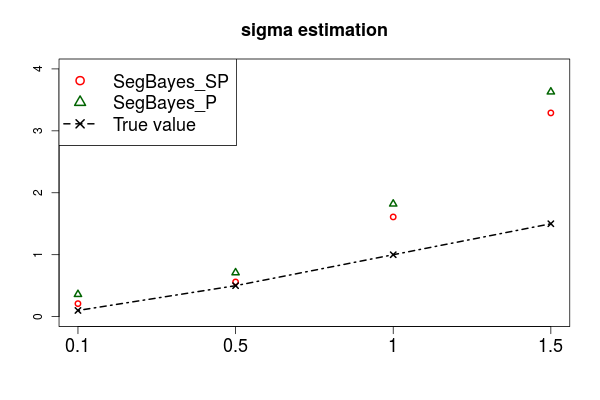}
    \caption{\small{Average of the estimates $\hat{\sigma}$ on 100 simulated series with $\sigma=0.1$, $\sigma=0.5$, $\sigma=1$ and $\sigma=1.5$, for {\it SegBayes\_SP} (semi-parametric model) and {\it SegBayes\_P} (parametric model).}}
    \label{simulations100graph_sigma}
    \end{center}
    \end{figure}

\paragraph{Overall results and comparison between \textit{SegBayes\_SP} and \textit{SegBayes\_P}.}
Figure \ref{simulations100graph} shows the average quality criteria over the 100 simulated series for the 4 different values of $\sigma$, and Figure \ref{simulations100graph_sigma} gives the average of the estimates $\hat{\sigma}$.

\begin{figure}[!htp]
    \begin{center}
    \includegraphics[scale=0.4]{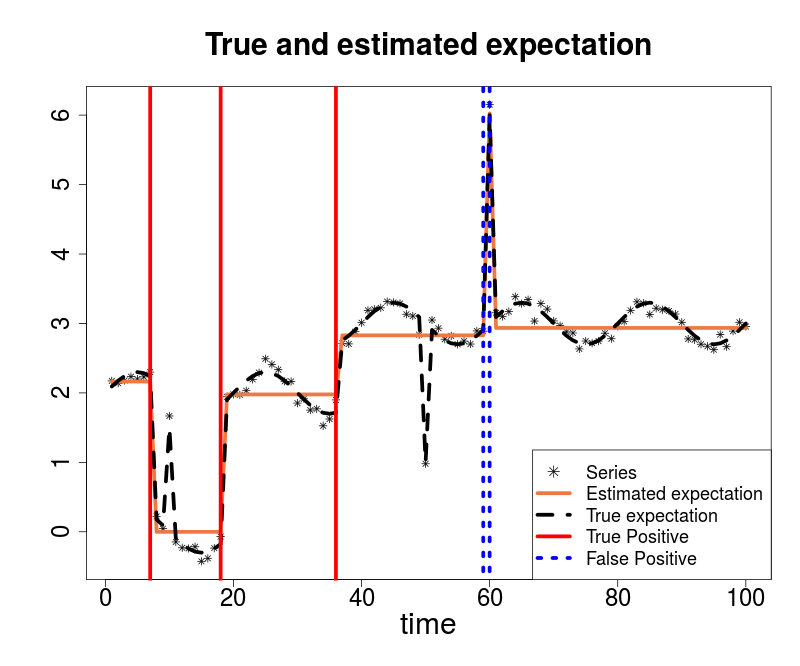}
    \caption{\small{Result of \textit{SegBayes\_P} on the particular simulated series with $\sigma=0.1$: the series, the whole true expectation and its estimation, the True Positive (TP), False Positive (FP) and False Negative (FN) change-points are also represented.}}
    \label{reconstruction_noise0.1_nobias}
    \end{center}
    \end{figure}

When the detection problem is easy ($\sigma=0.1$), \textit{SegBayes\_SP} tends to recover correctly the two parts of the model. Indeed the selected number of segments is close to the true one and the change-points are well positioned (small $FDRbp$ and $FNRbp$). The same occurs with the selected functions. When the functional part is not taken into account, using the procedure \textit{SegBayes\_P}, the segmentation needs to compensate: the number of segments is overestimated leading to bad estimation of the segmentation part (higher $RMSE(\mu)$ and $FDRbp$). This is illustrated in Figure \ref{reconstruction_noise0.1_nobias} on one simulated series where the procedure selects false positive change-points to fit well the series (in particular to catch the peaks of the function). The higher $FNRbp$ compared to \textit{SegBayes\_SP} is explained by the fact that the precise positioning will be better obtained when the functional part is well estimated, in particular at positions where the jump of the means is small.

When the noise $\sigma$ increases, the two procedures tend to underestimate the number of segments and for \textit{SegBayes\_SP} also the number of functions. These results were expected (and generally observed for the segments in segmentation problems). Indeed, in this case, the segmentation and the functional part will be confused together and also with the errors. To avoid false detection and selection, one may prefer to select less change-points and functions. This is particularly marked for the functional part for which we observed a big decreasing of the number of selected functions from $\sigma=0.5$. As a consequence, $RMSE(f)$ increases and $FNR(f)$ is almost close to $1$. Moreover since the functional is not recovered, the two procedures lead to the same results for the estimation of the segmentation part (same values of the different criteria for the quality of the segmentation).

Finally the quality of the estimation of both parts is related to the quality of the estimation of $\sigma$: estimates $\hat{\sigma}$ are good for $\sigma=0.1$ and $\sigma=0.5$, but they appear too large for $\sigma=1$ and $\sigma=1.5$ (see Figure \ref{simulations100graph_sigma}).

An advantage of these Bayesian procedures is that we obtain empirical posterior probabilities for the possible change-points and functions. 
This will be illustrated in the following paragraph.


\paragraph{Results for particular series.}

We look here in detail at particular series simulated as follows: three change-points are considered at positions 7, 18 and 36, the means over the four segments are 2, 0, 2 and 3 respectively and we consider two series: one with $\sigma=0.1$ and the second one with $\sigma=1$.

\begin{figure}[!htp]
      \begin{center}
      \includegraphics[scale=0.18]{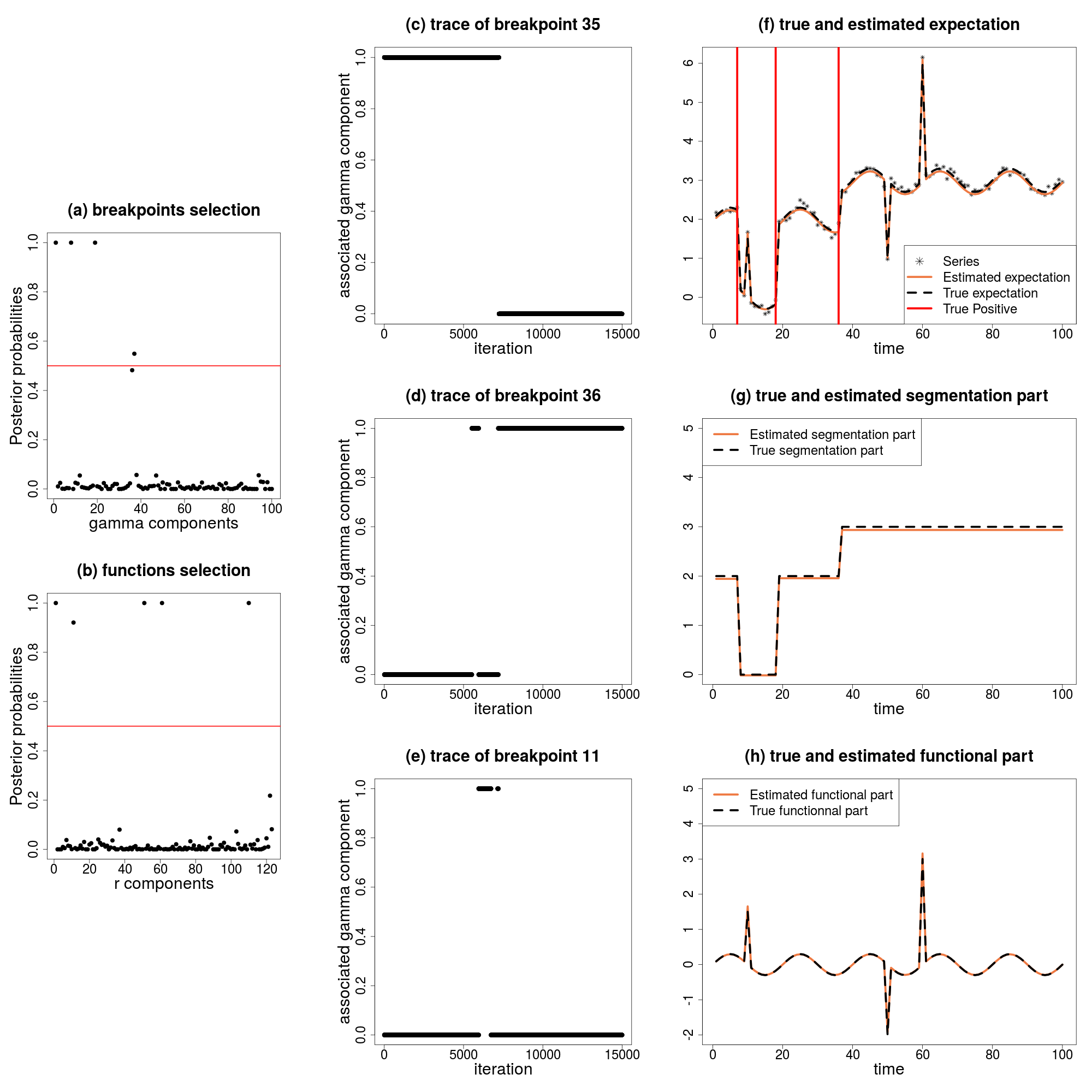}
      \caption{\small{Results of \textit{SegBayes\_SP} on the particular simulated series with $\sigma=0.1$. Posterior probabilities for the $\boldsymbol{\gamma}$ and $\mathbf{r}$ components (plots ($a$) and ($b$)). Traces of the $\boldsymbol{\gamma}$ components 35, 36 and 11 (plots ($c)$, ($d$) and ($e$)). The series, the whole true expectation and its estimation, the True Positive (TP), False Positive (FP) and False Negative (FN) change-points are also represented (plot ($f$)). The true and estimated segmentation part and functional part (plots ($g$) and ($h$)). }}
      \label{Exemple_SP_0.1}
      \end{center}
    \end{figure}

The results of the procedure \textit{SegBayes\_SP} for the first series with $\sigma=0.1$ (low noise) are given in Figure \ref{Exemple_SP_0.1}. We observe that, in this case, the true change-points and functions are exactly recovered (see plots ($a$) and ($b$)) and the two parts of the model are well estimated (see plots ($g$) and ($h$)) leading to a good estimation of the whole (see plot ($f$)).  Note that Table \ref{indicesFunctions} gives the correspondences between the selected functions from the dictionary and their indexes, and by convention the change-point at time 0 and the constant function are selected.\\
The advantage of a Bayesian approach compared to a frequentist one, is that the posterior probabilities give some additional interesting information. For instance, the posterior probabilities of the change-points at positions 0, 7 and 18 are 1, while those of the change-points 35 and 36 are 0.48 and 0.55 respectively (see plots ($a$) and ($b$)). Indeed, the choice between change-points 35 and 36 is not so easy for the sampler (see the traces in plots ($c$) and ($d$)). Observe also that we can deduce from the traces that the posterior probability of having a change-point in $\{35,36\}$ is close to 1.
Moreover, let us give another comment. Instead of using a threshold to obtain $\widehat{\boldsymbol{\gamma}}$ and $\widehat{\mathbf{r}}$, another way could have been to use the posterior modes. However, this approach generally leads to more false positive change-points. This has been observed in many series and this is also the case in this particular series. Here using the posterior modes  results in the detection of the true change-points at positions 7 and 18 as previously, but also of a change-point at position 35 instead of 36 and of a false positive change-point at position 11.

On the same simulated series, the result of the procedure \textit{SegBayes\_P} (without considering the functional part) is given in Figure \ref{reconstruction_noise0.1_nobias}. The change-points 7, 18, 36, 59 and 60 are selected with high posterior probabilities. The two last detected change-points are false positive but they correspond to the Haar 60 from the functional part. As explained in the previous paragraph, when the functional part is forgotten, the segmentation tends to catch it.

		\begin{table}[!htp]
		\begin{center}
		\begin{tabular}{|c|c|c|c|c|}
		\hline
		\cellcolor{lightgray} Index & 1 & $i\in\{2,\ldots,101\}$& $110$ & $120$ \\
		\hline
		\cellcolor{lightgray} Function & constant term & Haar function at $t=i-1$&$\sin\big(2\pi \times 5 \times \frac{t}{100}\big)$& $\sin\big(2\pi \times 10 \times \frac{t}{100}\big)$\\\hline
		\end{tabular}
		\end{center}
		\caption{\small{Functions selected from the dictionary and their corresponding indexes.}}
		\label{indicesFunctions}
		\end{table}

\begin{figure}[!htp]
      \begin{center}
    \includegraphics[scale=0.3]{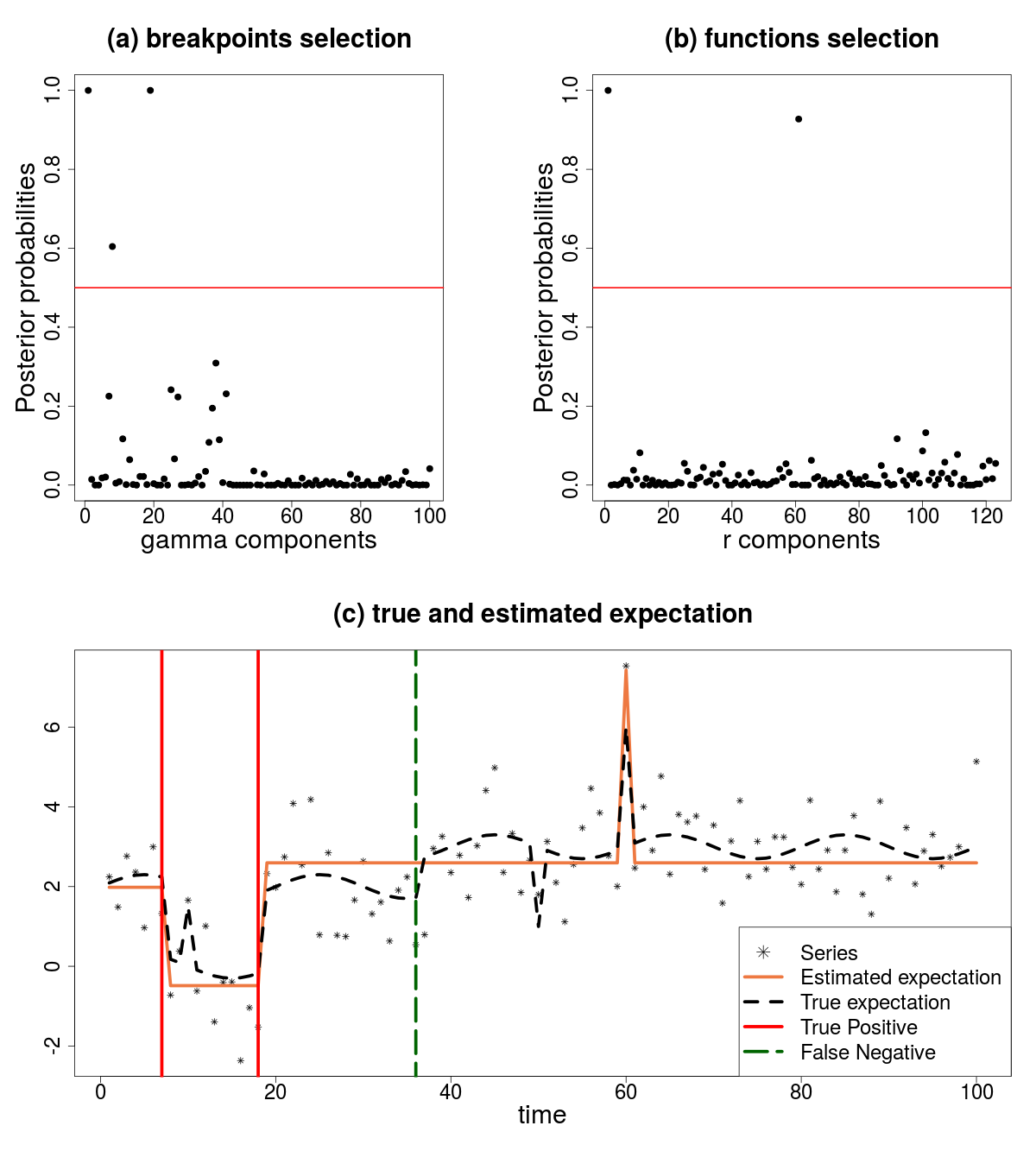}
      \caption{\small{Results of \textit{SegBayes\_SP} on the particular simulated series with $\sigma=1$. Top: posterior probabilities for the $\boldsymbol{\gamma}$ and $\mathbf{r}$ components. Bottom: the series, the whole true expectation and its estimation, the True Positive (TP), False Positive (FP) and False Negative (FN) change-points are also represented.}}
      \label{Example_SP_1}
      \end{center}
    \end{figure}

As pointed out previously, the case $\sigma=1$ is challenging since the jump of 1 (even 2) on the mean of the series become difficult to detect, as well as peaks functions which can be confounded with the noise. In Figure \ref{Example_SP_1}, we can see that using \textit{SegBayes\_SP} the posterior probabilities of the true change-points at positions 7 and 18 are 0.60 and 1 respectively, while the posterior probabilities of the change-point at position 36 which is not selected is 0.19. Concerning the functional part, only the Haar 60 is detected with posterior probability 1. Looking at the series, this is expected since at position 36, the jump is not marked. The result of the procedure \textit{SegBayes\_P} on this same series is not very good, since only the true change-point 18 is detected, with a posterior probability of 0.72.

\paragraph{Sensitivity and convergence.}

To study the sensitivity of the estimates $\widehat{\boldsymbol{\gamma}}$ and $\widehat{\mathbf{r}}$ to the choice of prior parameters, we ran the Metropolis-Hastings algorithm on the same particular series as before with $\sigma=0.1$, with different choices of prior parameters.
Table \ref{TableSensitivity} in Appendix \ref{App:AppendixC} shows the results obtained for several simulation scenarios (21 different runs).

On average, the procedure is not over-sensitive to the choice of the prior parameters: among the 21 runs, 10 detected exactly the true change-points and  functions, and 8 runs detected the true change-points and functions with a shift or an``exchange". For instance runs 9 and 13 select change-points at positions 37 and 35 respectively instead of the 36. Runs 6, 11, 16 and 19 select change-point at position 10 instead of 7, and functions 9 and 10 (Haar 8 and Haar 9) instead of 11 (Haar 10).

Some other sensitivity remarks can still be done. First we can see that too small values for $c_1$ and $c_2$ should not be used since it results on too many undetected change-points and functions (see run 2). Moreover the number of components of $\boldsymbol{\gamma}$ and $\mathbf{r}$ to be changed at each iteration should not be too high. Indeed, in this case, the proposed changes are too difficult to accept leading to a poor acceptance rate. For instance, for runs 12 and 15, the acceptance rates for $\boldsymbol{\gamma}$ and $\mathbf{r}$ are  respectively 0.05\% and 0.02\% (instead of 3\% for the other runs). It results that for these two runs some false-positive change-points and functions are selected. On the contrary, the initial number of segments and functions and the values of the probabilities $\pi_l$ and $\eta_j$ do not seem to influence too much the number and the validity of selected change-points and functions from the dictionary (see runs 4 to 9 and 16 to 21).

To study the convergence, we ran the algorithm three times with the same prior parameters than run 1, with 20000 iterations (5000 of burn-in), and one time with 50000 iterations (10000 of burn-in).
The results are given in Table \ref{TableConvergence} in Appendix \ref{App:AppendixC}. We observe that 20000 iterations including 5000 of burn-in seem enough to reach convergence since the obtained results are similar for these four runs (or to the previous runs 1, 5, 8, 11, 14, 17 and 20 that have the same prior parameters). The acceptance rates for $\boldsymbol{\gamma}$ and $\mathbf{r}$ are not very high in general (around 3\%). But if we look in more detail at the traces, it appears that usually when a true change-point or a true function from the dictionary is selected, it will be selected until the end of the algorithm, while a position which is not a change-point will be alternatively selected and unselected. As a consequence, when most of the change-points and bias functions have been selected, the chain will not be much updated, resulting in a poor acceptance rate.

\section{Application}\label{sec:appli}

\subsection{Geodetic data}

In this section, we propose to use our procedure in the geodesic field for the problem of homogenization of GPS series. Indeed, such series are used to determine accurate station velocities for tectonic and Earths mantle studies \citep{King2010}. However they are affected by two effects: (i) abrupt changes that are related to instrumental changes (documented or not), earthquakes or changes in the raw data processing strategy and (ii) periodic signals that are due to environmental signals, such as soil moisture or atmospheric pressure changes. The correct detection of these effects is fundamental for the aforementioned application.

Here we consider a particular series (the height coordinate of the series) from the GPS station in Yarragadee, Australia, YAR2 at the weekly scale. The data can be downloaded at
{\footnotesize{\verb+http://sideshow.jpl.nasa.gov/post/series.html+.}}
We refer the reader to \cite{BCLM:2014} for more details about the problem and the data.

We apply our proposed procedure {\it{SegBayes\_SP}} to this series with a dictionary of $194$ functions that includes the constant function and the Fourier functions:
$t\mapsto\sin\left(2\pi w_i t\right),t\mapsto\cos\left(2\pi w_i
t\right)$ where $w_i=i/T$, $T=\text{max}(t)-\text{min}(t)$ and $T/i$
is larger than 8 weeks.
The Metropolis-Hastings algorithm is run for 100 000 iterations (30 000 burn-in), with $c_1=c_2=50$. The initial number of segments and functions is 5, the number of change-points or functions proposed to be changed at each iteration is 1. The initial probability for each possible function is 0.01, the initial probability for a position not associated with a known equipment change or malfunction is 0.01 and the initial probability for a position associated with a known equipment change or malfunction is 0.5.
Concerning the Gibbs sampler algorithm, we run it for 100 000 iterations including 50000 burn-in iterations and we choose $c_1=c_2=50$.

\begin{figure}[!h]
	\begin{center}
	\includegraphics[width=12cm,height=6cm]{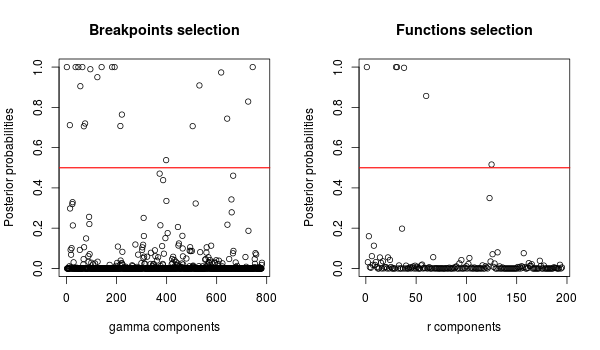}
	\includegraphics[scale=0.2]{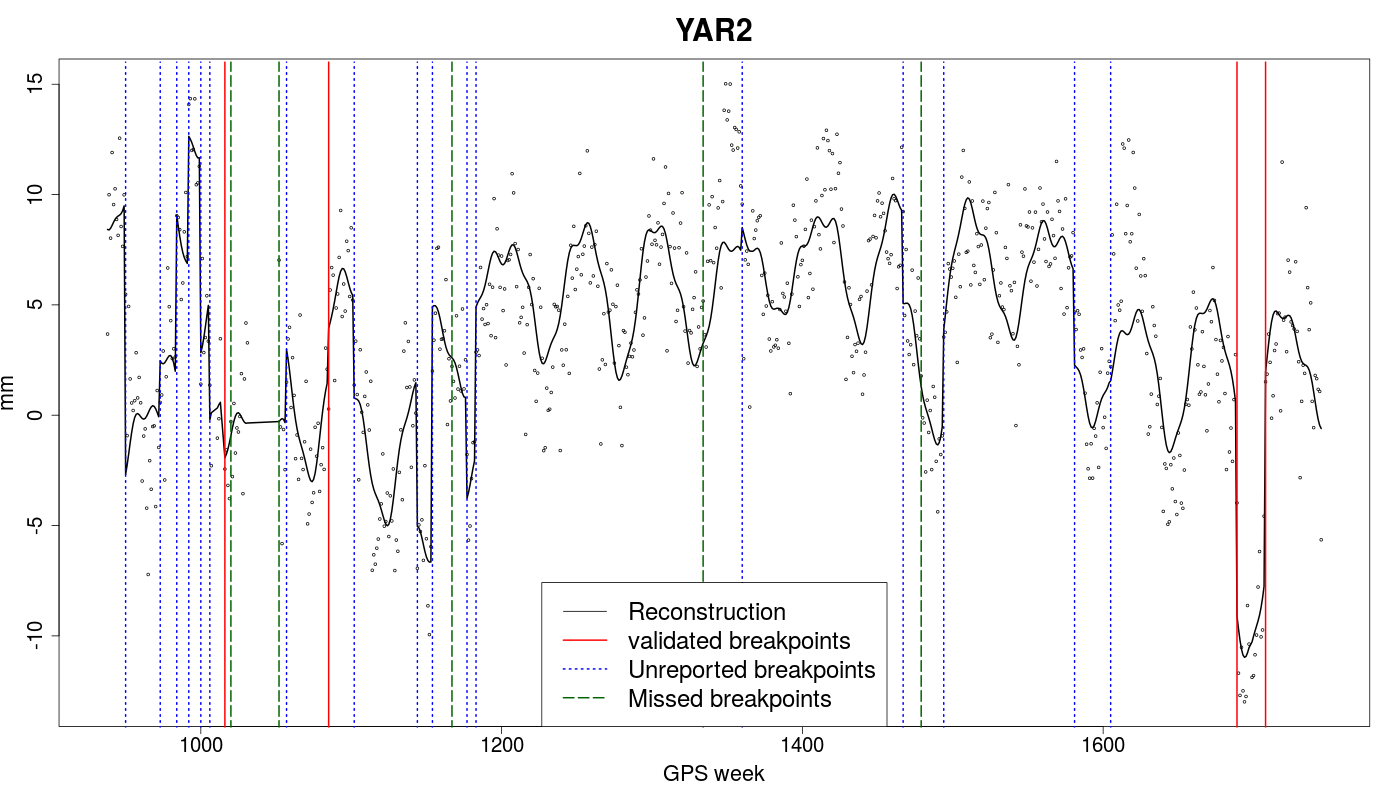}
	\caption{\small{Result of the procedure \textit{SegBayes\_SP} on the YAR2 series. Top: posterior probabilities of the change-points and functions. Bottom: Estimated expectation and validated, unreported and missed change-points (based on known equipment changes and malfunctions).}}
\label{Results_YAR2}
	\end{center}
	\end{figure}

In the top of Figure \ref{Results_YAR2}, the posterior probabilities of the change-points and functions are given. In the bottom of Figure \ref{Results_YAR2}, the series, its reconstruction,  the validated, unreported and missed change-points are plotted (validated means reported in databases and detected, unreported means not reported in databases but detected, and missed means reported in databases but non detected). The different status of the change-points should be tempered since they are based on known equipment changes and malfunctions. Some unknown malfunctions should have appeared, hence some unreported change-points could in fact be real change-points. Moreover small earthquakes can have been not reported in databases.

Four change-points previously reported in databases are detected by \textit{SegBayes\_SP}. In particular, GPS weeks 1016 and 1085 correspond exactly to clock changes, and GPS weeks 1689 and 1708 correspond exactly to radom radar changes. The unreported change-point at GPS week 1057 maybe associated with the receiver and clock change at GPS week 1031. Note that we also apply our procedure with non informative prior for the positions of the change-points (all the initial probability equal to 0.01). In this case, several documented breakpoints are not yet detected, showing the gain of using prior knowledge.

Examining the selected functions, 4 of them were selected by the procedure: $\sin\big(2\pi \times \frac{t}{52}\big)$,  $\cos\big(2\pi \times \frac{t}{52}\big)$,  $\sin\big(2\pi \times \frac{t}{41.05}\big)$ and $\sin\big(2\pi \times \frac{t}{26}\big)$. These functions furnish relevant geodetic information. In particular, the selection of the three functions with periods of 52 and 26 weeks is consistent with the fact that atmospheric pressure can be approximated by periodic signals with dominant annual and semi-annual periods (\citealp{DFCM2002}).

\subsection{Exchange rate data }

In this section, we propose to use our procedure in the econometrics field for the problem of the exchange rate. More precisely, we study the daily records of Mexican Peso/US Dollar exchange rate from January 2007 to December 2012 (data available at $\verb+ www.federalreserve.gov+$). These data were studied by \cite{Martinez:14} with a Bayesian nonparametric method.

We apply our proposed procedure {\it{SegBayes\_SP}} to this series with a dictionary of $23$ functions that includes the following functions of time $\{t\mapsto t^j,\quad j=1,2\}$
and the Fourier functions:
$t\mapsto\sin\left(2\pi i t/n \right),t\mapsto\cos\left(2\pi i
t/n \right)$ for $i=1,\ldots,10$. The Metropolis-Hastings algorithm is run for 100 000 iterations (30 000 burn-in), with $c_1=c_2=50$. The initial number of segments and functions are 5, the number of change-points or functions proposed to be changed at each iteration is 1.

The initial probability for each possible function is 0.01, as well as the initial probability for a position.
Concerning the Gibbs sampler algorithm, we run it for 100 000 iterations including 50000 burn-in iterations and we choose $c_1=c_2=50$.

\begin{figure}[!h]
	\begin{center}
	\includegraphics[width=12cm,height=6cm]{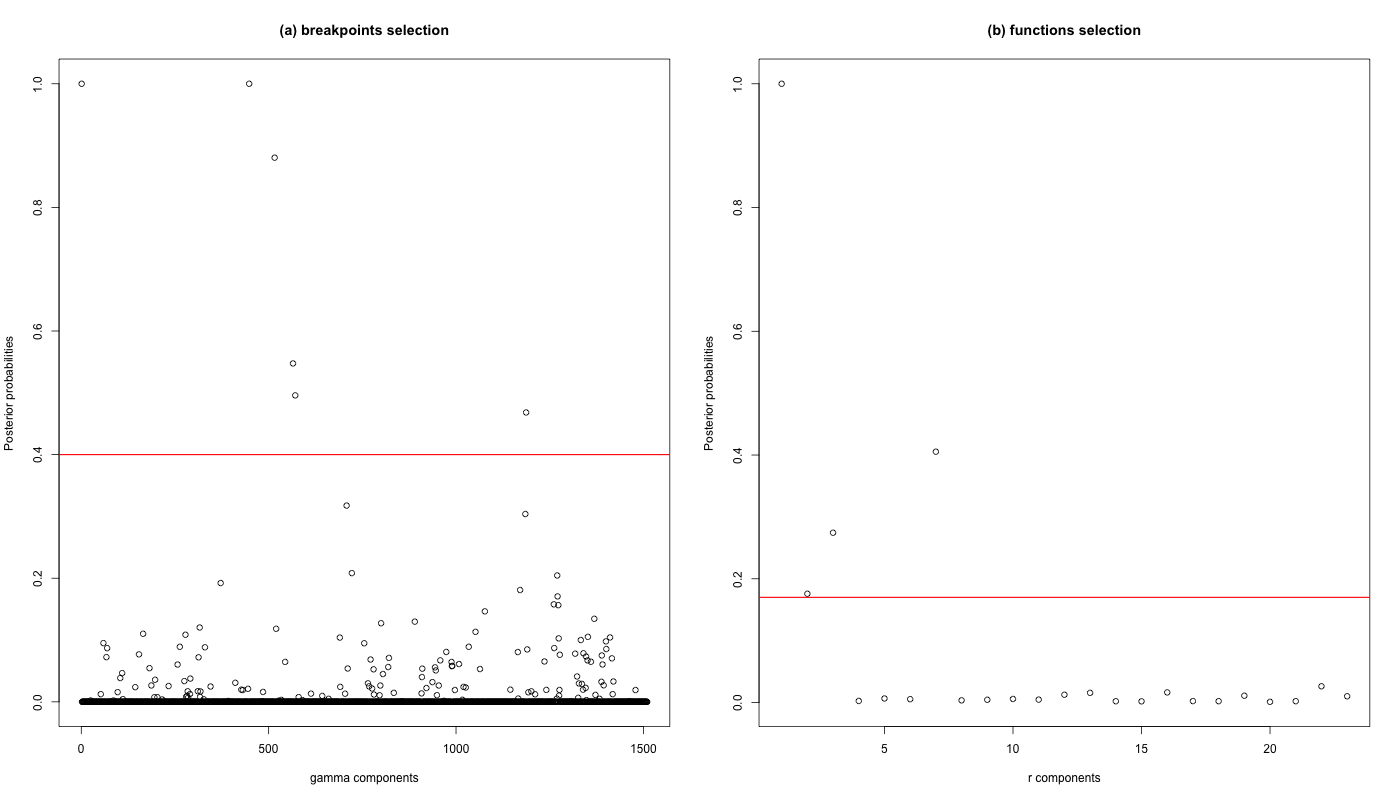}
	\includegraphics[scale=0.25]{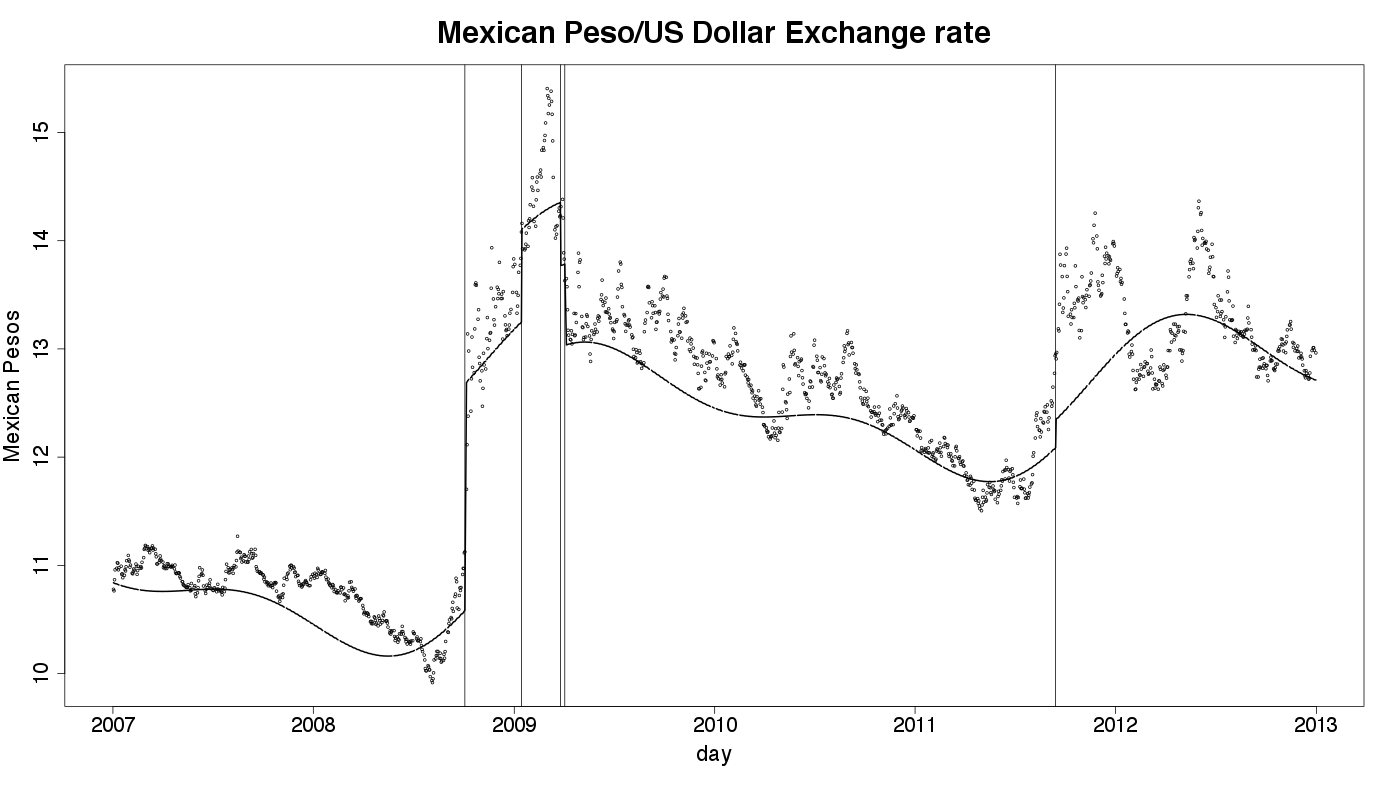}
	\caption{\small{Result of the procedure \textit{SegBayes\_SP} on the Mexican Peso/US Dollar exchange rate series. Top: posterior probabilities of the change-points and functions. Bottom: Estimated expectation and change-points.}}
\label{Results_Exchange}
	\end{center}
	\end{figure}

In the top of Figure \ref{Results_Exchange}, the posterior probabilities of the change-points and functions are given. In the bottom of Figure \ref{Results_Exchange}, the series with both the estimation of the expectation and the estimated change-points. The five most probable change-points detected are at dates October 03, 2008, January 14, 2009, March 26, 2009,  April 3, 2009 and September 14, 2011.

These are close to the ones obtained by \cite{Martinez:14} and, as explained in this latter paper, four can be related to events in Mexico and USA: 1) September-October 2008: pick of the 2007-2008 financial crisis; 2) March-May 2009: the flu pandemic suffered by Mexico; 3) the US debt-ceiling crisis in 2011. The other change-points detected by \cite{Martinez:14}, that do not correspond to known events, are not detected by our procedure, certainly due to the presence of the functional part.

\section{Discussion} \label{sec:discussion}

In this paper, we propose a novel Bayesian method to segment one-dimensional piecewise-constant signals corrupted by a functional part. The functional part is estimated by a linear combination of functions from a dictionary. Since the dictionary can be large, this method is flexible and allows us to estimate functions with both smooth components and local irregularities. For the estimation of the change-points, following \cite{HL:2010}, the piecewise-constant term is expressed as a product of a lower triangular matrix by a sparse vector. This presents the advantage that the Stochastic Search Variable Selection approach can be applied, avoiding the use of the Reversible Jump method in which the mixing of the chains can be problematic (see \citealp{LL2000}).\\
The global estimation procedure is based on a Metropolis-Hastings algorithm for estimating the positions of the change-points and for selecting the functions (estimation of the latent vectors $\boldsymbol{\gamma}$ and $\mathbf{r}$), 
and on a Gibbs sampler to estimate the parameters $\boldsymbol{\beta}_{\boldsymbol{\gamma}}$, $\boldsymbol{\lambda}_{\mathbf{r}}$ and $\sigma^2$, conditionally on $\boldsymbol{\gamma}$ and $\mathbf{r}$. It is quite common in Stochastic Search Bayesian Variable Selection to integrate all the parameters except the latent vectors, see for instance \cite{GMC1997} or \cite{BotRi2010}. We obtain good results in this paper in a simulation study and on two real data sets. In particular, we illustrate in the simulation study that taking into account the functional part improves the detection of the breakpoints.

Our procedure benefits from the Bayesian framework, which results on two important aspects. The first one is that posterior distributions of the parameters are obtained. Then, we can have a quantification of the uncertainty through posterior probabilities for the possible change-points and functions composing the functional part, or through credible intervals. Credible intervals have not been used  in our examples, but they can be easily calculated for the means over the segments $(\mu_1,\ldots,\mu_K)$, for the coefficients $(\lambda_1,\ldots,\lambda_M$) of the selected functions or for the functional part $f$.
The second important aspect is that we can introduce expert knowledge in the model through prior distributions (see Section \ref{sec:appli}). In particular in the geodetic field, some change-points are not detected when non informative priors are used whereas they are detected when previous knowledge is taking into account. Moreover, change-points and functions composing the functional part are detected in an automatic way. This could provide a significant improvement for the users in this field since previously the detection of change-points was done visually.

An important issue is the choice of the criterion to estimate the parameters $\boldsymbol{\gamma}$ and $\mathbf{r}$. 
As proposed by \cite{MPRR04}, the criterion used in this work minimizes a loss function which is the sum of the False Discovery and of the False Negative, leading to a threshold of $1/2$ for the posterior probabilities.  In the simulation study, this criterion seems to outperform other strategy based on the posterior mode which leads to more false-positive change-points. Moreover other thresholds can be used, to minimize different loss functions, see \cite{MPRR04} or \cite{MPR06}. Finally, an other way to select the final change-points and functions from the dictionary would be to run the algorithm say three times, and to take the intersections of the three results, for both the change-points and the bias functions. That would lead to perfect results for most of the groups of three runs from the sensitivity analysis.

To use our procedure, hyper-parameters should be chosen, but the sensitivity analysis shows that the procedure is not over-sensitive to these choices.

Eventually, we think that this kind of approach should not be restricted to the problem of detecting change-points in a one-dimensional series contaminated by a signal. Indeed, it should be interesting and useful to extend this approach to the joint segmentation of multiple series corrupted by a functional part.

\section*{Acknowledgments}
Meili Baragatti, Karine Bertin, Emilie Lebarbier and Cristian Meza are supported by FONDECYT grants 1141256 and 11412578, and the mathamsud 16-MATH-03 SIDRE project.

\appendix
\section{Integration of the joint posterior distribution} \label{App:AppendixA}
Integrating the joint posterior (\ref{jointposterior}) with respect to $\boldsymbol{\beta}_{\boldsymbol{\gamma}}$ we obtain:
\begin{eqnarray*}
&\pi(\boldsymbol{\gamma},\boldsymbol{\lambda}_{\mathbf{r}},\mathbf{r},\sigma^2|\mathbf{Y})&\propto  (2\pi\sigma^2)^{-n/2}(1+c_1)^{-d_{\boldsymbol{\gamma}}/2}\pi(\boldsymbol{\gamma})\pi(r)\sigma^{-2}\\
&&\times \exp\left[ -\frac{1}{2\sigma^2}(\mathbf{Y} -\mathbf{F}_{\mathbf{r}}\boldsymbol{\lambda}_{\mathbf{r}})^\prime  \left\{ I-\frac{c_1}{1+c_1} \mathbf{X}_{\boldsymbol{\gamma}} (\mathbf{X}_{\boldsymbol{\gamma}}^\prime \mathbf{X}_{\boldsymbol{\gamma}})^{-1} \mathbf{X}_{\boldsymbol{\gamma}}^\prime \right\} (\mathbf{Y}-\mathbf{F}_{\mathbf{r}}\boldsymbol{\lambda}_{\mathbf{r}}) \right]\\
&&\times (2\pi)^{-d_{\mathbf{r}}/2} \left| c_2 \sigma^2 (\mathbf{F}_{\mathbf{r}}^\prime \mathbf{F}_{\mathbf{r}})^{-1} \right|^{-1/2}\exp \left[-\frac{1}{2\sigma^2} \boldsymbol{\lambda}_{\mathbf{r}}^\prime \left(\frac{\mathbf{F}_{\mathbf{r}}^\prime \mathbf{F}_{\mathbf{r}}}{c_2} \boldsymbol{\lambda}_{\mathbf{r}}\right) \right].
\end{eqnarray*}
Integrating with respect to $\boldsymbol{\lambda}_{\mathbf{r}}$, we obtain:
\begin{eqnarray*}
\pi(\boldsymbol{\gamma},\mathbf{r},\sigma^2|\mathbf{Y})&\propto & (2\pi\sigma^2)^{-\frac{n}{2}} (1+c_1)^{-d_{\boldsymbol{\gamma}}/2}\pi(\boldsymbol{\gamma})\pi(r)\sigma^{-2} \left(\frac{\left| \left( \mathbf{F}_{\mathbf{r}}^\prime \left(\mathbf{U}_{\boldsymbol{\gamma}}^{-1}+\frac{I}{c_2} \right) \mathbf{F}_{\mathbf{r}}\right)^{-1}  \right|}{|c_2(\mathbf{F}_{\mathbf{r}}^\prime \mathbf{F}_{\mathbf{r}})^{-1}|} \right)^{1/2}\\
&&\times \exp\left[ -\frac{1}{2\sigma^2} \mathbf{Y}^\prime\left( \mathbf{U}_{\boldsymbol{\gamma}}^{-1}-\mathbf{U}_{\boldsymbol{\gamma}}^{-1} \mathbf{F}_{\mathbf{r}}\left(\mathbf{F}_{\mathbf{r}}^\prime \left(\mathbf{U}_{\boldsymbol{\gamma}}^{-1}+\frac{I}{c_2} \right) \mathbf{F}_{\mathbf{r}}\right)^{-1}\mathbf{F}_{\mathbf{r}}^\prime \mathbf{U}_{\boldsymbol{\gamma}}^{-1} \right) \mathbf{Y}  \right].
\end{eqnarray*}
Finally, integrating over $\sigma^2$, we obtain the integrated posterior (\ref{integratedposterior}).

\section{Tables of the sensitivity study} \label{App:AppendixC}

\begin{table}
\tiny{
		\rotatebox{90}{\begin{tabular}{|c|c|c|c|c|c|c|c|c|c|}
		\hline
		\rowcolor{lightgray}  & \small{Values} & \small{Initial} & \small{Initial} & \small{Nb of $\boldsymbol{\gamma}$}  & \small{Nb of $\mathbf{r}$} & \small{Values} & \small{Values} & \small{Selected} & \small{Selected}\\	
		\rowcolor{lightgray} \small{Run} & \small{of $c_1$ } & \small{nb of} & \small{nb of} & \small{comp.} & \small{comp.} & \small{of the} & \small{of the} & \small{change-points} & \small{functions}\\
		\rowcolor{lightgray}  & \small{and $c_2$} & \small{segments} & \small{functions} & \small{changed} & \small{changed} & \small{$\pi_l$} & \small{$\eta_j$} & & \small{(function of index 1 always}\\
		\rowcolor{lightgray}  &  & &  & \small{at each} & \small{at each } &  & & & \small{selected (constant term))}\\
	\rowcolor{lightgray} 	&&&& \small{ iter.} &  \small{ iter.} &&&&\\
		\hline
		1 & 50 &  &  &  &  &  &  & \footnotesize{7, 18, 36} & \footnotesize{11, 51, 61, 110}\\
		2 & 10 &  &  &  &  &  &  & \footnotesize{7, 18} & \footnotesize{61}\\
		3 & 500 &  \multirow{-3}{*}{3} & \multirow{-3}{*}{3} & \multirow{-3}{*}{2} & \multirow{-3}{*}{2} & \multirow{-3}{*}{1/100} & \multirow{-3}{*}{1/100} & \footnotesize{7, 18, 36}  & \footnotesize{11,51,61,110}\\
		\hline
		4 &  & 1 &  &  &  &  &  & \footnotesize{7, 18, 36, 49, 50}  & \footnotesize{11, 61, 110}\\
		5 &  & 3 &  &  &  &  &  & \footnotesize{7, 18, 36}  & \footnotesize{11, 51, 61, 110}\\
		6 & \multirow{-3}{*}{50} &  10 & \multirow{-3}{*}{3} & \multirow{-3}{*}{2} & \multirow{-3}{*}{2} & \multirow{-3}{*}{1/100} & \multirow{-3}{*}{1/100} & \footnotesize{10, 18, 36} & \footnotesize{9, 10, 51, 61, 110}\\
		\hline
		7 &  &  & 1 &  &  &  &  &  \footnotesize\footnotesize{7, 18, 36} & \footnotesize{11, 51, 61, 110}\\
		8 &  &  & 3 &  &  &  &  & \footnotesize{7, 18, 36} & \footnotesize{11, 51, 61, 110}\\
		9 & \multirow{-3}{*}{50} & \multirow{-3}{*}{3} & 10 & \multirow{-3}{*}{3} & \multirow{-3}{*}{2} & \multirow{-3}{*}{1/100} & \multirow{-3}{*}{1/100} & \footnotesize{7, 18, 37}& \footnotesize{11, 51, 61, 110}\\
		\hline
		10 &  &  &  & 1 &  &  &  & \footnotesize{7, 18, 36} & \footnotesize{11, 51, 61, 110}\\
		11 &  &  &  & 2 &  &  &  & \footnotesize{10, 18, 36} & \footnotesize{9, 10, 51, 61, 110}\\
		12 & \multirow{-3}{*}{50} & \multirow{-3}{*}{3} & \multirow{-3}{*}{3} & 5 & \multirow{-3}{*}{2} & \multirow{-3}{*}{1/100} & \multirow{-3}{*}{1/100} & \footnotesize{7, 18, 35, 74, 82, 98} & \footnotesize{11, 37, 51, 61, 110}\\
		\hline
		13 &  &  &  &  & 1 &  &  & \footnotesize{7, 18, 35} & \footnotesize{11, 51, 61, 110}\\
		14 &  &  &  &  & 2 &  &  & \footnotesize{7, 18, 36} & \footnotesize{11, 51, 61, 110}\\			
		15 & \multirow{-3}{*}{50} & \multirow{-3}{*}{3} & \multirow{-3}{*}{3} & \multirow{-3}{*}{2} & 5 & \multirow{-3}{*}{1/100} & \multirow{-3}{*}{1/100} & \footnotesize{7, 18, 36, 59, 60} & \footnotesize{11, 16, 51, 68, 110, 120}\\	
		\hline
		16 &  &  &  &  &  & 1/20 &  & \footnotesize{10, 18, 36} & \footnotesize{9, 10, 51, 61, 110}\\
		17 &  &  &  &  &  & 1/100 &  & \footnotesize{7, 18, 36} & \footnotesize{11, 51, 61, 110}\\
		18 & \multirow{-3}{*}{50} & \multirow{-3}{*}{3} & \multirow{-3}{*}{3} & \multirow{-3}{*}{2} & \multirow{-3}{*}{2} & 1/500 & \multirow{-3}{*}{1/100} & \footnotesize{7, 18, 36} & \footnotesize{11, 51, 61, 110}\\	
		\hline
		19 &  &  &  &  &  &  & 1/20 & \footnotesize{10, 18, 36} & \footnotesize{9, 10, 51, 61, 110}\\
		20 &  &  &  &  &  &  & 1/100 & \footnotesize{7, 18, 36} & \footnotesize{11, 51, 61, 110}\\
		21 & \multirow{-3}{*}{50} & \multirow{-3}{*}{3} & \multirow{-3}{*}{3} & \multirow{-3}{*}{2} & \multirow{-3}{*}{2} & \multirow{-3}{*}{1/100} & 1/500 & \footnotesize{7, 18, 36, 59, 60} & \footnotesize{11, 51, 110}\\
		\hline
		\end{tabular}}
		\caption{\small{\textit{SegBayes\_SP}: prior parameters used in the different runs of the Metropolis-Hastings algorithm applied on the particular series with $\sigma=0.1$. The number of iterations is 20000 with a burn-in of 5000 iterations for all runs.}}
		}
		\label{TableSensitivity}
		\end{table}

\begin{table}[h!]
		\begin{center}
		\begin{tabular}{|c|c|c|c|c|}
		\hline
		\rowcolor{lightgray} & Number of & burn-in & Selected & Selected \\
		\rowcolor{lightgray} Run & iterations &  & change-points & functions\\
		\hline
		22 & 20000 & 5000 & 7, 18, 36 & 1, 11, 51, 61, 110\\
		23 &20000 & 5000 & 7, 18, 36 & 1, 51, 61, 76, 110\\
		24 &20000 & 5000 & 8, 18, 36 & 1, 9, 11, 51, 61, 110\\
		25 &50000 & 10000 & 7, 18, 36 & 1, 11, 51, 61, 110\\
		\hline
		\end{tabular}
		\end{center}
		\caption{\small{\textit{SegBayes\_SP}: results of four runs of the Metropolis-Hastings algorithm applied on the particular series with $\sigma=0.1$, with the same prior parameters than run 1 in Table 2.}}
		\label{TableConvergence}
	    \end{table}

\newpage

\begin{thebibliography}{0}
\providecommand{\natexlab}[1]{#1}
\providecommand{\url}[1]{\texttt{#1}}
\expandafter\ifx\csname urlstyle\endcsname\relax
  \providecommand{\doi}[1]{doi: #1}\else
  \providecommand{\doi}{doi: \begingroup \urlstyle{rm}\Url}\fi

\end{thebibliography}


\begin{thebibliography}{99}

\newcommand{\enquote}[1]{``#1''}
\expandafter\ifx\csname natexlab\endcsname\relax\def\natexlab#1{#1}\fi
\expandafter\ifx\csname url\endcsname\relax
  \def\url#1{{\tt #1}}\fi
\expandafter\ifx\csname urlprefix\endcsname\relax\def\urlprefix{URL }\fi
\ifx\endbibitem\undefined \let\endbibitem\relax\fi


\bibitem[{Barbieri and Berger(2004)}]{BB2004}
Barbieri, M.M. and Berge, J.O. (2004).
\newblock \enquote{Optimal predictive model selection.}
\newblock {\em The Annals of Statistics\/}, 32(3): 870--897.
\endbibitem

\bibitem[{Bertin et~al.(2016)Bertin, Collilieux, Lebarbier, and
  Meza}]{BCLM:2014}
Bertin, K., Collilieux, X., Lebarbier, E., and Meza, C. (2016).
\newblock \enquote{{Segmentation of multiple series using a Lasso strategy}.}
\newblock Submitted - arXiv:1406.6627.
\endbibitem

\bibitem[{Bickel et~al.(2009)Bickel, Ritov, and Tsybakov}]{tsybakov}
Bickel, P.~J., Ritov, Y., and Tsybakov, A.~B. (2009).
\newblock \enquote{Simultaneous analysis of lasso and {D}antzig selector.}
\newblock {\em Ann. Statist.\/}, 37(4): 1705--1732.
\endbibitem

\bibitem[{Bottolo and Richardson(2010)}]{BotRi2010}
Bottolo, L. and Richardson, S. (2010).
\newblock \enquote{Evolutionary stochastic search for Bayesian model exploration.}
\newblock {\em Bayesian Analysis\/}, 5(3): 583--618.
\endbibitem

\bibitem[{Boys and Henderson(2004)}]{BH:2004}
Boys, R.~J. and Henderson, D.~A. (2004).
\newblock \enquote{{A Bayesian approach to DNA sequence segmentation}.}
\newblock {\em Biometrics\/}, 60: 573--588.
\newblock With discussion.
\endbibitem

\bibitem[{Caussinus and Mestre(2004)}]{CM2004}
Caussinus, H. and Mestre, O. (2004).
\newblock \enquote{{Detection and correction of artificial shifts in climate
  series}.}
\newblock {\em Applied Statistics\/}, 53: 405--425.
\endbibitem


\bibitem[{Dobigeon et~al.(2007)Dobigeon, Tourneret, and Scargle}]{DTS:2007}
Dobigeon, N., Tourneret, J.-Y., and Scargle, J. (2007).
\newblock \enquote{{Joint segmentation of multivariate astronomical time
  series: bayesian sampling with a hierarchical model}.}
\newblock {\em IEEE Transactions on Signal Processing\/}, 55: 414--423.
\endbibitem

\bibitem[{Dong et~al.(2002)Dong, {P. Fang}, Bock, Cheng, and
  Miyazaki}]{DFCM2002}
Dong, D., {P. Fang}, P., Bock, Y., Cheng, M.~K., and Miyazaki, S. (2002).
\newblock \enquote{{Anatomy of apparent seasonal variations from GPS-derived
  site position time series}.}
\newblock {\em Journal of Geophysical Research (Solid Earth)\/}, 107(B4): ETG
  9--1--ETG 9--16.
\endbibitem

\bibitem[{Fearnhead(2006)}]{F:2006}
Fearnhead, P. (2006).
\newblock \enquote{{Exact and efficient Bayesian inference for multiple
  changepoint problems}.}
\newblock {\em Statistics and Computing\/}, 16: 203--213.
\endbibitem

\bibitem[{George and McCulloch(1993)}]{GeorgeMcCulloch}
George, E. and McCulloch, R. (1993).
\newblock \enquote{Variable selection via Gibbs sampling.}
\newblock {\em Journal of the American Statistical Association\/}, 88(423):
  881--889.
\endbibitem

\bibitem[{Georges and McCulloch(1997)}]{GMC1997}
Georges, E.I. and McCulloch, R.E. (1997).
\newblock \enquote{Approaches for Bayesian variable selection.}
\newblock {\em Statistica Sinica\/}, 7: 339--373.
\endbibitem
%

\bibitem[{Harchaoui and L\'evy-Leduc(2010)}]{HL:2010}
Harchaoui, Z. and L\'evy-Leduc, C. (2010).
\newblock \enquote{Multiple Change-Point Estimation With a Total Variation
  Penalty.}
\newblock {\em Journal of the American Statistical Association\/}, 105:
  1480--1493.
\endbibitem




\bibitem[{H{\"a}rdle et~al.(1998)H{\"a}rdle, Kerkyacharian, Picard, and
  Tsybakov}]{picard}
H{\"a}rdle, W., Kerkyacharian, G., Picard, D., and Tsybakov, A. (1998).
\newblock {\em Wavelets, approximation, and statistical applications\/}, volume
  129 of {\em Lecture Notes in Statistics\/}.
\newblock Springer-Verlag, New York.
\endbibitem

\bibitem[{King et~al.(2010)King, Altamimi, Boehm, Bos, Dach, Elosegui, Fund,
  Hernandez-Pajares, Lavall{\'e}e, Cerveira, Riva, Steigenberger, van Dam,
  Vittuari, Williams, and Willis}]{King2010}
King, M., Altamimi, Z., Boehm, J., Bos, M., Dach, R., Elosegui, P., Fund, F.,
  Hernandez-Pajares, M., Lavall{\'e}e, D., Cerveira, P., Riva, R.,
  Steigenberger, P., van Dam, T., Vittuari, L., Williams, S., and Willis, P.
  (2010).
\newblock \enquote{{Improved constraints on models of glacial isostatic
  adjustment}.}
\newblock {\em A review of the contribution of ground-based geodetic
  observations, Surveys in Geophysics\/}, 31(5): 465--507.
\endbibitem


\bibitem[{Lavielle and Lebarbier(2001)}]{LL2000}
Lavielle, M. and Lebarbier, E. (2001).
\newblock \enquote{An application of MCMC methods for the multiple change-points problem.}
\newblock {\em Signal processing \/}, 81: 39--53.
\endbibitem


\bibitem[{Liu(1994)}]{Liu}
Liu, J.S (1994).
\newblock \enquote{The collapsed Gibbs sampler in Bayesian computations with application to a gene regulation problem.}
\newblock {\em J. Am. Stat. Ass.\/}, 89(427): 958--966.
\endbibitem


\bibitem[{Mart{\'i}nez, and Mena(2014)}]{Martinez:14}
Mart{\'i}nez, A. and Mena, R.(2014).
\newblock \enquote{On a Nonparametric Change Point Detection
Model in Markovian Regimes.}
\newblock {\em Bayesian Analysis\/}, 9: 823--858.
\endbibitem


\bibitem[{Muller, Parmigiani, Robert and Rousseau(2004)}]{MPRR04}
Muller, P., Parmigiani, G., Robert, C., and Rousseau, J.(2004).
\newblock \enquote{Optimal sample size for multiple testing: the case of gene expression microarrays.}
\newblock {\em J. Amer. Statist. Assoc.\/}, 99: 990--1001.
\endbibitem

\bibitem[{Muller, Parmigiani and Rice(2006)}]{MPR06}
Muller, P., Parmigiani, G. and Rice, K. (2006).
\newblock \enquote{FDR and Bayesian multiple comparisons rules.}
\newblock {\em Proc Valencia / ISBA 8th World meeting on Bayesian Statistics.\/}.
\endbibitem

\bibitem[{Picard et~al.(2011)Picard, Lebarbier, Budinska, and Robin}]{PLBR11}
Picard, F., Lebarbier, E., Budinska, E., and Robin (2011).
\newblock \enquote{{Joint segmentation of multivariate Gaussian processes using
  mixed linear models}.}
\newblock {\em Comp. Stat. and Data Analysis\/}, 55: 1160--1170.
\endbibitem


\bibitem[{Ruggieri(2013)}]{R:2013}
Ruggieri, E. (2013).
\newblock \enquote{{A Bayesian approach to detecting change points in climatic
  records}.}
\newblock {\em International Journal of Climatology\/}, 33: 520--528.
\endbibitem

\bibitem[{Smith and Kohn(1997)}]{SmithKohn}
Smith, M. and Kohn, R. (1997).
\newblock \enquote{{Non parametric regression using Bayesian variable
  selection}.}
\newblock {\em Journal of Econometrics\/}, 75: 317--344.
\endbibitem

\bibitem[{Tai and Xing(2010)}]{Tai:2010}
Tai, T.L. and Xing, H. (2010).
\newblock \enquote{A fast Bayesian change point analysis for the segmentation of microarray data.}
\newblock {\em Bioinformatics.\/}, 24: 2143--2148.
\endbibitem

\bibitem[{van Dyk and Park (2008)}]{vanDyk}
{van Dyk}, D.A. and Park, T. (2008).
\newblock \enquote{Partially collapsed Gibbs samplers: theory and methods.}
\newblock {\em J. Am. Stat. Ass.\/}, 103: 790--796.
\endbibitem

\bibitem[{Williams(2003)}]{W2003}
Williams, S. (2003).
\newblock \enquote{{Offsets in Global Positioning System time series}.}
\newblock {\em Journal of Geophysical Research (Solid Earth)\/}, 108(B6, 2310).
\endbibitem

\bibitem[{Wyse et~al.(2011)Wyse, Friel, and Rue}]{WFR:2011}
Wyse, J., Friel, N., and Rue, H. (2011).
\newblock \enquote{{Approximate simulation-free Bayesian inference for multiple
  changepoint models with dependence within segments}.}
\newblock {\em Bayesian Analysis\/}, 6: 501--528.
\endbibitem


\bibitem[{Zellner(1986)}]{ZE1986}
Zellner, A. (1986).
\newblock \enquote{Bayesian inference and decision techniques -- essays in honour of Bruno De Finetti, chapter On assessing prior distributions and Bayesian regression analysis with $g$-prior distributions.}
\newblock {\em Goel, P.K. and Zellner, A\/}, 233--243.
\endbibitem


\end{thebibliography}

\end{document}